\documentclass[a4paper,12pt]{article}
\usepackage[cp1251]{inputenc}
\usepackage[english, russian, ukrainian]{babel}
\usepackage{amsmath}
\usepackage{amssymb}

\newcommand{\1}{1\!\!\,{\rm I}}

\newcommand{\vf}{\varphi}

\newcommand{\mbR}{{\mathbb R}}
\newcommand{\mfX}{{\mathfrak X}}

\renewcommand{\tilde}{\widetilde}

\newcommand{\pt}{\partial}

\begin{document}
 \large
UDC 519.21 \vskip10pt

\begin{center}
{\bf
Smoothing problem in anticipating scenario}\\[1cm]
Andrey A. Dorogovtsev\\[1 cm]
Institute of Mathematics, Ukrainian Academy of Sciences,\\
ul. Tereshenkovskaia,3,Kiev, Ukraine\\
e-mail: adoro@imath.kiev.ua
\end{center}
\vskip1cm

{\bf Introduction.} This article is devoted to the stochastic
anticipating equations with the extended stochastic integral with
respect to the Gaussian processes of a special type and its
application to the smoothing problem in the case when noise is
represented by the two jointly Gaussian Wiener processes, which can
have not a semimartingale property with respect to the joint
filtration. In order to describe the objects of our consideration
more explicitly consider the following example.

{\bf Example 0.1}. Consider the ordinary stochastic differential
equation in $\mbR$
$$
dx(t)=a(x(t))dt+b(x(t))dw(t)
$$
with the smooth enough coefficients $a$ and $b.$  Denote by $x(r,s,t)$ the
solution which starts at the moment $s$ from the point $r.$ Let the function
$\vf\in C^2(\mbR)$  has bounded derivatives. Define for $r\in\mbR,
s\in[0;T]$
$$
\Phi(r,s)=\Gamma(\vf(r,s,T)),\eqno(0.1)
$$
where $\Gamma$ is the certain operator of the second quantization [2]. In
particular $\Gamma$  can be a mathematical expectation. Then, it can be
proved, that $\Phi$  satisfies the following partial stochastic differential
equation
$$
d\Phi(s,r)=-\left[
\frac{1}{2}b^2(r)\frac{\pt^2}{\pt r^2}\Phi(r,s)+a(r)\frac{\pt}{\pt r}\Phi(r,
s)\right] ds+
$$
$$
+b(r)\frac{\pt}{\pt r}\Phi(r,s)d\gamma(s).
$$
Here $\gamma(s)=\Gamma w(s)$  and the last differential is treated
in the sense of anticipating stochastic integration. When $\Gamma$
is the mathematical expectation, the last term vanishes.

This example shows the main goal of this article. Namely there exist
situations when the naturally arising Wiener functionals satisfy the
anticipating stochastic differential equations and can be described
with the using of stochastic calculus. Here we propose the
appropriate machinery and derive the correspondent equations.

We will consider the second quantization transformation of the
different Wiener functionals. For such transformed functionals we
will get the anticipating stochastic equations with the extended
stochastic integral. Accordingly to this aim the article is
organized as follows. The first part contains the properties of the
second quantization operators in connection with the extended
stochastic integral or, more generally, with the Gaussian strong
random operators [3].  The section 2 and 3 are devoted to the
following pair of equations
$$
\begin{aligned}
dx_1(t)=a_1(x_1(t))+dw_1(t)\\
dx_2(t)=a_2(x_1(t))+dw_2(t),
\end{aligned}
$$
where $w_1, w_2$ are jointly Gaussian Wiener processes, which can
have not a semimartingale property with respect to the joint
filtration. Here we will look for the equation for
$E(f(x_1(t))/x_2).$

 {\bf 1. Second quantization and integrators.}   The material
of this section is partially based on the works [4,5]. Corresponding
facts are placed here for the completeness of the exposition but
their proofs are omitted. New claims are presented with the proofs.

We will start here with the abstract picture, when the ``white
noise'' generated by the Wiener process is substituted by the
generalized Gaussian random element in the Hilbert space. Let $H$
be a separable real Hilbert space with the   norm $\|\cdot\|$ and
inner product $(\cdot,\cdot).$  Suppose that $\xi$ is the
generalized Gaussian random element in $H$ with zero mean and
identical covariation. In other words $\xi$ is the family of
jointly Gaussian random variables denoted by $(\vf, \xi), \vf\in
H$ with the properties

1) $(\vf, \xi)$ has the normal distribution with zero mean and
variance $\|\vf\|^2$  for every $\vf\in H,$

2) $(\vf, \xi)$  is linear with respect to $\vf.$

During this section we suppose that all the random variables and
elements are measurable with respect to $\sigma(\xi)=\sigma((\vf,
\xi), \vf\in H.$ If the random variable $\alpha$ has the finite
second moment, than $\alpha$ has an Ito-Wiener expansion [6]
$$
\alpha=\sum^\infty_{k=0}A_k(\xi, \ldots,\xi). \eqno(1.1)
$$
Here, for every $k\geq1 \ A_k(\xi, \ldots, \xi)$   is the infinite-dimensional
generalization of the Hermite polinomial from $\xi,$   correspondent to the
$k$-linear symmetric Hilbert-Shmidt form $A_k$  on $H.$  Moreover, now the
following relation holds
$$
E\alpha^2=\sum^\infty_{k=0}k!\|A_k\|^2_k.
\eqno(1.2)
$$
Here $\|\cdot\|_k$  is the Hilbert-Shmidt form in $H^{\otimes k}.$  The same
expansion for $H$-valued random elements will be necessary. Let $x$  be a
random element in $H$  such, that
$$
E\|x\|^2<+\infty.
$$
Then, for every $\vf\in H$
$$
(x,\vf)=\sum^\infty_{k=0}A_k(\vf; \xi,\ldots,\xi).
\eqno(1.3)
$$
It can be easily checked using (1.2), that now $A_k$  is the
$k+1$-linear (not necessary symmetric) Hilbert-Shmidt form. So one
can write now
$$
x=\sum^\infty_{k=0}\tilde{A}_k(\xi,\ldots,\xi), \eqno(1.4)
$$
where
$$
\tilde{A}_k(\vf_1,\ldots,\vf_k):=
\sum^\infty_{j=1}A_k(e_j; \vf_1,\ldots, \vf_k)e_j
$$
for the arbitrary orthonormal basis $\{e_j; j\geq1\}$  in $H$   and the
series (1.4)  converges in $H$ in the square mean. The relation (1.2)
remains to be true
$$
E\|x\|^2=\sum^\infty_{k=0}k!\|\tilde{A}_k\|^2_k,
\eqno(1.5)
$$
where $\|\tilde{A}_k\|_k$  is the Hilbert--Shmidt norm in the space of
$H$-valued $k$-linear forms on $H.$

Now recall the definition of the operators of the second
quantization. Let $C$  be a continuous linear operator in $H.$
Suppose that the operator norm $\|C\|\leq1.$   Then for $\alpha$
and $x$  from (1.1)  and (1.4) define
$$
\begin{aligned}
&\Gamma(C)\alpha =\sum^\infty_{k=0}A_k(C\xi,\ldots,C\xi),\\
&\Gamma(C)x =\sum^\infty_{k=0}\tilde{A}_k(C\xi,\ldots,C\xi),\\
\end{aligned}
\eqno(1.6)
$$
where for $k\geq1 \ A_k(C\cdot,C\cdot,\ldots,C\cdot)$ and
$\tilde{A}_k(C\cdot,C\cdot,\ldots,C\cdot)$ are new Hilbert-Shmidt forms.

Using the estimation
$$
\|A_k(C\cdot,C\cdot,\ldots,C\cdot)\|_k\leq\|C\|^k\cdot\|A_k\|_k
$$
it is  easy to prove [2], that $\Gamma(C)$ is a continuous linear
operator in the space of square integrable random variables or
elements in $H.$

{\bf Definition 1.1.} [2] Operator $\Gamma(C)$  is the operator of \
t\,h\,e \ s\,e\,c\,o\,n\,d \
q\,u\-\,a\,n\,\-t\,i\,\-z\,a\,\-t\,i\,o\,n \ correspondent to the
operator $C.$

 Before to consider some examples, we will present the useful
representation of the  second quantization operators. Let $\xi'$  be
the generalized Gaussian random element in $H,$  independent and
equidistributed with $\xi.$

Consider the following generalized Gaussian random element in $H$
$$
\eta=\sqrt{1-CC^*}\xi'+C\xi.
\eqno(1.7)
$$
This element can be properly defined by the formula
$$
\forall \vf\in H: \ \ \ \ \
(\vf, \eta):=(\sqrt{1-CC^*}\vf,\xi')+(C^*\vf,\xi).
\eqno(1.8)
$$
Note that $\eta$  has zero mean and identity covariation. In order to check
this,  is sufficient to note the relation
$$
\|\vf\|^2=\|\sqrt{1-CC^*}\vf\|^2+\|C^*\vf\|^2.
$$
For every random variable $\alpha$ with an expansion (1.1)  define
$$
\alpha(\eta):=\sum^\infty_{k=0}A_k(\eta,\ldots,\eta).
$$
The following representation will be useful.

{\bf Lemma 1.1.} {\sl For arbitrary $\alpha\in L_2(\Omega, \sigma(\xi), P)$
and operator $C$ in $H$ with $\|C\|\leq1$
$$
\Gamma(C)\alpha =E(\alpha(\eta)/\xi).
\eqno(1.9)
$$
}

{\it Proof}. Note, that the both parts of (1.9) are continuous
with respect to  $\alpha$ in the square mean. So, it is enough to check
(1.9)  for the following random variables
$$
e^{(\vf,\xi)-\frac{1}{2}\|\vf\|^2}, \ \vf\in H. \eqno(1.10)
$$
Really,  the random variable of this kind has the Ito-Wiener
expansion of the form
$$
e^{(\vf,\xi)-\frac{1}{2}\|\vf\|^2}=\sum^\infty_{k=0}\frac{1}{k!}
\vf^{\otimes k}(\xi,\ldots,\xi).
$$
Here $\vf^{\otimes k}$ is $k$-th tensor power of $\vf$  which acts on
$H$ by the rule
$$
\vf^{\otimes k}(\psi_1,\ldots,\psi_k)=\prod^k_{j=1}(\vf,\psi_j).
$$
So, as it was mentioned in introduction, for $\alpha,$ which has an
expansion (1.1)
$$
E\alpha e^{(\vf,\xi)-\frac{1}{2}\|\vf\|^2}=\sum^\infty_{k=0}A_k
(\vf,\ldots,\vf).
$$
Hence, the set of all linear combinations of the variables (1.10)
is dense in $L_2.$ Now [2] the following equality holds
$$
\Gamma(C)e^{(\vf,\xi)-\frac{1}{2}\|\vf\|^2}=
e^{(C^*\vf,\xi)-\frac{1}{2}\|C^*\vf\|^2}. \eqno(1.11)
$$
From other side
$$
e^{(\vf,\eta)-\frac{1}{2}\|\vf\|^2}=
e^{(C^*\vf,\xi)-\frac{1}{2}\|C^*\vf\|^2}\cdot
e^{(\sqrt{1-CC^*}\vf,\xi')-\frac{1}{2}\|\sqrt{1-CC^*}\vf\|^2}.
$$
In order to finish the proof it is enough now note, that
$$
E
e^{(\sqrt{1-CC^*}\vf,\xi')-\frac{1}{2}\|\sqrt{1-CC^*}\vf\|^2}=1,
$$
and that $\xi'$ and $\xi$  are independent. It follows from here, that
$$
E\left(e^{(\vf,\eta)-\frac{1}{2}\|\vf\|^2}/\xi\right)=
e^{(C^*\vf,\xi)-\frac{1}{2}\|C^*\vf\|^2}.
$$
Lemma is proved.

This lemma has the following useful application for us.

{\bf Corollary 1.1.} {\sl
Let $\Gamma(C)$ be an operator of the second quantization. Let $x$  be $a$
random element in the complete separable metric space $\mfX$
measurable with respect to $\xi.$   Then there exists the random
probability measure $\mu$ on $\mfX$  such, that for every bounded
measurable function $f: \mfX\to\mbR$   the following equality holds
$$
\int_\mfX fd\mu =\Gamma(C)f(x).
$$
}

{\it Proof.} Let us define $\mu$  as a conditional distribution of $x(\eta)$
with respect to $\xi.$  Then, for every bounded measurable $f: \mfX\to\mbR$
$$
\int_\mfX fd\mu=E(f(x(\eta))/\xi)=\Gamma(C)f(x).
$$
The unique difficulty on this way lies in the proper definition of
$x(\eta)$ (remind, that $\xi$ and $\eta$ are not usual random
elements). In order to break this difficulty we will use the
following analog of the Levy theorem. Let $\{e_j; j\geq1\}$   be
an orthonormal basis in $H.$  Define the sequences of random
elements in $H$  by the rule
$$
\begin{aligned}
&
\xi_n=\sum^n_{j=1}(e_j,\xi)e_j,\\
&
\eta_n=\sum^n_{j=1}(e_j,\eta)e_j, \ n\geq1.
\end{aligned}
$$
Note that the sequences $\{\xi_n;n\geq1\}$  and $\{\eta_n;n\geq1\}$ are
equidistributed. Now for every $n\geq1$ consider the random measure
$\nu_n$ in $\mfX,$  which is built in the following way
$$
\nu_n(\Delta)=E\{\1_\Delta(x)/\xi_n\}.
$$
Here $\Delta$ is an arbitrary Borel subset of $\mfX.$  This random
measures have two important properties. First of all for every
$n\geq1$ $\nu_n$ can be viewed as $\widetilde{\nu}_n(\xi_n),$
where $\widetilde{\nu}_n$ is a Borel function from $H$ to the
space of all probability measures on $\mfX$ equipped with the
distance of weak convergence. Secondly, with probability one
$\nu_n$ weakly converge to $\delta_x$ under $n$ tends to infinity.
The last assertion follows from the usual Levy theorem [7]. More
precisely, for arbitrary continuous bounded function $f:
\mfX\to\mbR$
$$
f(x)=E(f(x)/\xi)=\lim_{n\to\infty}E(f(x)/\xi_n)=
$$
$$
=\lim_{n\to\infty}\int_\mfX f(u)\nu_n(du) \ \mbox{a.s.}
$$
Taking $f$  from the countable set which define the weak
convergence [8] we get the required statement. Now note, that the
sequence of random measures $\{\widetilde{\nu}_n(\eta\; n\geq1\}$
is equidistributed with $\{\widetilde{\nu}_n(\xi_n); n\geq1\}.$
Hence with probability one there exists the weak limit of
$\widetilde{\nu}_n(\eta_n)$  which is a delta-measure concentrated
in the certain random point $y.$  This random point $y$ is by
definition $x(\eta).$  The correctness of this definition can be
easily checked. Lemma is proved.

Consider the examples of the random measures, which arise in the
application of the Corollary 1.1 and will be important for us.

{\bf Example 1.1.}  Suppose, that $H=L_2([0;T], \mbR^d).$  Define the
generalized Gaussian random element $\xi$ in $H$  with the help of the
$d$-dimensional Wiener process $W$ on $[0; T].$  Namely, for
$
\vf=(\vf_1,\ldots,\vf_d)\in L_2([0;1], \mbR^d)$ define
$$
(\vf,\xi):=\sum^d_{j=1}\int^T_0\vf_j(s)dW_j(s).
\eqno(1.12)
$$
Now consider the functions $a: \mbR^d\to\mbR^d$ and $b:
\mbR^d\to\mbR^{d\times d}$ which satisfy the Lipshitz condition
and the domain $G$ in $\mbR^d$  with the $C^1$-boundary $\Gamma$.
Let for every $s\in[0;T]$  and $u\in\mbR^d$ $x(u,s,T)$  denote the
solution at time $T$ of the following Cauchy problem
$$
\begin{cases}
dx(t)=a(x(t))dt+b(x(t))dW(t),\\
x(s)=u.
\end{cases}
\eqno(1.13)
$$
Denote by $\nu_{u,s}$  the random measure obtained from $x(u,s,T)$
via corollary 1.1 with the help of the certain operator of the
second quantization $\Gamma(C).$  In the next section we will obtain
the stochastic variant of the Kolmogorov equation for $\nu_{u,s}.$
Note, that in the case $C=0$ measures $\nu_{u,s}$ became to be
deterministic and satisfy the usual Kolmogorov equation [1].

Now let us define for every $u\in G$  the random moment
$$
\tau_{u,s}=\inf\{T, t\leq T: x(u,s,t)\in\Gamma\}.
$$
Let $\mu_{u,s}$ be the random measure obtained from
$x(u,s,\tau_{u,s})$ via corollary 1.1. It occurs, that measures
$\mu_{u,s}$ satisfy certain anticipating boundary  value problem.

In order to describe the anticipating SPDE for the random measures
from above mentioned example we need in the relation between the
operators of the second quantization and extended stochastic
integral. We will study this connection in the more general
situation when the extended stochastic integral is substituted by
the general Gaussian strong random operator (GSRO in the sequel).
Let us recall the following definition.

{\bf Definition 1.2}. [3] T\,h\,e \ G\,a\,u\,s\,s\,i\,a\,n \
s\,t\,r\,o\,n\,g \ r\,a\,n\,d\,o\,m \ l\,i\,n\,e\,a\,r \
o\,p\,e\,\-r\,a\,\-t\,o\,r \ (\,G\,S\,R\,O\,) \ $A$ in $H$ is the
mapping, which maps every element $x$ of $H$  into the jointly
Gaussian with $\xi$ random element in $H$ and is continuous in the
square mean.

As an example of GSRO the integral with respect to Wiener process can be
considered.

{\bf Example 1.2.} Consider $H$  and $\xi$  from example 1.1. Let for
simplicity $d=1.$  Define GSRO $A$ in the following way
$$
\forall \vf\in H: \ \ \ \ \ \
(A\vf)(t)=\int^t_0\vf(s)dw(s), \ \ t\in[0; T].
$$
It can be easily seen that $A\vf$ now is a Gaussian random element in
$H,$  and $A$ is continuous in square mean.

For to include in this picture the  integration with respect to
another Gaussian processes (for example with respect to the
fractional Brownian motion) consider more general GSRO. Suppose,
that $K$ be a bounded linear operator, which acts from
$L_2([0;T])$  to $L_2([0;T]^2).$   Define
$$
\forall \vf\in H: \ \ \ \
(A\vf)(t)=\int^T_0(K\vf)(t,s)dw(s).
$$
It can be checked, that $A$ is GSRO in $H.$  Making an obvious changes one can
define the GSRO acting from the different Hilbert space $H_1$  into $H.$
For example consider for $\alpha\in\left(
\frac{1}{2};1\right)$ the covariation function of the fractional Brownian
motion [9]  with Hurst parameter $\alpha$
$$
R(s,t)=\frac{1}{2}(t^{2\alpha}+s^{2\alpha}-|t-s|^{2\alpha}).
$$
Define the space $H_1$  as a completion of the set of step
functions on $[0;T]$ with respect the inner product under which
$$
(\1_{[0;s]},\1_{[0;t]})=R(s,t).
$$
Consider the kernel $K^\alpha$ from the integral representation of the
fractional
Brownian motion $B^\alpha$ [10]
$$
B^\alpha(t)=\int^t_0 K^\alpha(t,s)dw(s)
$$
and
$$
\frac{\pt K^\alpha}{\pt t}(t,s)=c_\alpha\left(\alpha-\frac{1}{2}\right)
(t-s)^\alpha-\frac{3}{2}\left(\frac{s}{t}\right)^{\frac{1}{2}-\alpha}.
$$
Define for $\vf\in H_1$
$$
(K\vf)(t,s)=\int^t_s\vf(r)\frac{\pt K^\alpha}{\pt r}(r,s)dr\1_{[0;t]}(s).
$$
Now let
$$
(A\vf)(t)=\int^T_0(K\vf)(t,s)dw(s)=\int^t_0(K\vf)(t,s)dw(s).
$$
Then
$$
(A\vf)(t)=\int^t_0\vf(s)dB^\alpha(s).
$$

We will consider the action of GSRO on the random elements in $H.$
Corresponding definition was proposed in [3, 11]. Consider
arbitrary GSRO $A$ in $H.$ Then for every $\vf\in H$ the
Ito-Wiener expansion of $A\vf$ contains only two terms
$$
A\vf=\alpha_0\vf+\alpha_1(\vf)(\xi).\eqno(1.14)
$$
Here $\alpha_0$  is a continuous linear operator in $H$  and $\alpha_1$ is
a continuous linear operator from $H$ to the space of Hilbert-Shmidt
operators in $H.$ Now let $x$  be a random element in $H$  with the finite
second moment. Then $\alpha_1(x)$ has a finite second moment in the space
of Hilbert-Shmidt operators. So for every $\vf\in H$
$$
\alpha_1(x)(\vf)=\sum^\infty_{k=0}B_k(\vf; \xi,\ldots,\xi).
$$
It can be easily verified, that $B_k$ is $k+1$-linear $H$-valued
Hilbert-Shmidt form on $H.$ Define $\Lambda B_k$ as a symmetrization of
$B_k$  with respect to all $k+1$ variables.

{\bf Definition 1.3.} [3, 11]  The random element $x$ \
b\,e\,l\,o\,n\,g\,s \ t\,o \ t\,h\,e \ d\,o\,m\,a\,i\,n \ o\,f \
d\,e\,f\,i\,n\,i\,t\,i\,o\,n \ o\,f \ G\,S\,R\,O \ $A$ if the series
$$
\sum^\infty_{k=0}\Lambda B_k(\xi,\ldots,\xi)
$$
converges in $H$  in the square mean and in this case
$$
Ax=\alpha_0x+\sum^\infty_{k=0}\Lambda B_k(\xi,\ldots,\xi).
\eqno(1.15)
$$

In the partial cases this definition gives us the definition of
the extended stochastic integral [6, 12, 13]. We will define this
integral for the special class of Gaussian processes. Suppose,
that $H=L_2([0;T])$  and $\xi$ is generated by the Wiener process
$w$ as above. Consider the jointly Gaussian with $w$  process
$\{\gamma(t); t\in[0;T]\}$ with zero mean.

{\bf Definition 1.4.} [14]  Process $\gamma$ is \ a\,n \
i\,n\,t\,e\,g\,r\,a\,t\,o\,r \ if there exists the constant $C$ such
that for every step function $\vf$ on $[0;T]$
$$
\vf=\sum^{n-1}_{k=0}a_k\1_{[t_k; t_{k+1})}
$$
the following unequality holds
$$
E(\sum^{n-1}_{k=0}a_k(\gamma(t_{k+1})-\gamma(t_k)))^2\leq C
\sum^{n-1}_{k=0}a_k^2(t_{k+1}-t_k).
\eqno(1.16)
$$

The good examples of the integrators can be obtained via the following
simple  statement.

{\bf Lemma 1.2.} {\sl Let $\Gamma(C)$ be an operator of the second
quantization. Then
$$
\gamma(t)=\Gamma(C)w(t), t\in[0;T]
$$
is an integrator.
}

The proof of this lemma easily follows from the properties of $\Gamma(C).$

Note that the integrator can have the unbounded quadratic
variation and consequently can have not the semimartingale
properties (see [14]). It is easy to see from (1.16), that for
every integrator $\gamma$ and $\vf\in L_2([0;T])$ the stochastic
integral
$$
\int^t_0\vf d\gamma
$$
exists as a limit of the integrals from the step functions and
$$
E\left(\int^t_0\vf d\gamma\right)^2\leq C\int^t_0\vf^2(s)ds.
$$
So one can define GSRO $A_\gamma$  associated with the integrator $\gamma$
by the rule
$$
\forall \vf\in L_2([0;T]): \ \ \ \ \
(A_\gamma\vf)(t)=\int^t_0\vf d\gamma.
$$
In this situation the definition 1.3 became to be a definition of the
extended stochastic integral with respect to $\gamma.$  Note that in the
case $\gamma=w$ it will be a usual extended integral.

Now let us consider the relation between the action of GSRO and the
operators of the second quantization.

{\bf Theorem 1.1.} [4]  {\sl  Let $A$ be a GSRO in $H$ and
$\Gamma(C)$ be an operator of the second quantization. Suppose
that the random element $x$  lies in the domain of definition of
$A$ in the sence of definition 1.3. Then $\Gamma(C)x$ belongs to
the domain of definition of GSRO $\Gamma(C)A$  and the following
equality holds
$$
\Gamma(C)(Ax)=\Gamma(C)A (\Gamma(C)x).
\eqno(1.17)
$$
Here $\Gamma(C)A$ is the GSRO which acts by the rule
$$
\forall\vf\in H: \ \Gamma(C)A\vf=\Gamma(C)(A\vf).
$$
}

The proof of this theorem is placed in [4] and so is omitted.
Instead the proof consider the following important example of
application of the theorem 1.1 to the stochastic integration.

{\bf Example 1.3.}  Consider in the situation of the example 1.2
GSRO of integration with respect to Wiener process $w.$   Suppose
that random function $x$ in $L_2([0;T])$  with the finite second
moment is adapted to the flow of $\sigma$-fields generated by $w.$
It is well-known [12, 13], that in this case the extended
stochastic integral
$$
\int^t_0x(s)dw(s), \ t\in[0;T]
$$
exists and coincides with the  Ito integral. Now the theorem 1.1
says us that
$$
\Gamma(C)\left(\int^t_0x(s)dw(s)\right)=\int^t_0\Gamma(C)x(s)d\gamma(s),
$$
where $\gamma$ is an integrator of the type
$\gamma(t)=\Gamma(C)w(t)$  and the integral in the right part is
an extended stochastic integral.

{\bf2. Smoothing problem.}  The last sections of the article are
devoted to the following problem. Let $(w_1,w_2)$ be the pair of
jointly Gaussian one-dimensional Wiener processes. Let the processes
$x_1, x_2$ are obtained via the relations
$$
\begin{array}{l}
dx_1(t)=a_1(x_1(t))dt+dw_1(t),\\
dx_2(t)=a_2(x_1(t))dt+dw_2(t),\\
x_1(0)=x_2(0)=0.
\end{array}
\eqno(2.1)
$$

Note, that the second equality is just a definition of $x_2$ but not
an equation. The problem is to find the conditional distribution of
$x_1(t)$ for $t\in[0;1]$ under given $\{x_2(s); s\in[0;1]\}.$
 We will try to get the equation for
 $$
 E(f(x_1(t))/x_2)
 $$
 for the appropriate functions $f.$

 Firstly let us study the joint distribution of $(w_1,w_2)$. Note,
 that there exists the bounded linear operator $V: L_2([0;1])\to
 L_2([0;1])$ such, that \newline
 $\forall\vf_1,\vf_2\in L_2([0;1]):$
 $$
 E\int^1_0\vf_1dw_1\int^1_0\vf_2dw_2=\int^1_0\vf_1 V\vf_2 ds.
 $$

 This fact follows from the reason, that the left part of the above
 formula is the continuous bilinear form with respect to $\vf_1$ and
 $\vf_2.$  Moreover, the operator norm $\|V\|\leq 1.$

  In
this section we consider the density of the distribution $(x_1,
x_2)$  with respect to the distribution of $(w_1,w_2)$ and study its
properties under the conditional expectation. The problem is that
the distribution of $(w_1,w_2)$ is not a Wiener measure in $C([0;1],
\mbR^2).$  So in order to get the density we need to adapt the
general Gaussian measure setup [15] to the our case. For  the future
let us denote $C([0;1])$ as $C$  and identify the space $C([0;1],
\mbR^2)$  with the direct sum $C\oplus C,$ which is furnished by the
sum of the norms. Denote also by $H$  the space
$$
L_2([0;1],\mbR^2)=L_2([0;1])\oplus L_2([0;1])
$$
with the scalar product defined by the formula
$$
(\vf, \psi)=\int^1_0\vf_1\psi_1ds+\int^1_0\vf_2\psi_2ds.
$$
With the pair $(w_1,w_2)$  we can associate the generalized Gaussian
random element $\xi$ in $H$  by the rule
$$
(\vf, \xi)=\int^1_0\vf_1dw_1+\int^1_0\vf_2dw_2.
$$
Note that $\xi$  has not an identity covariation operator. Really
$$
E(\vf,\xi)(\psi, \xi)=\int^1_0\vf_1\psi_1ds+\int^1_0\vf_2\psi_2ds+
$$
$$
+\int^1_0\vf_1 V\psi_2ds +\int^1_0\psi_1 V\vf_2ds.
$$
Here $V$  is described above bounded linear operator in
$L_2([0;1]).$ Denote by $S$ the operator in $H$  which acts by the
rule
$$
S\vf=(\vf_1+V\vf_2, V^*\vf_1+\vf_2).
$$
Then
$$
E(\vf, \xi)(\psi,\xi)=(S\vf, \psi).
$$
Our aim is to describe the transformations of the pair $(w_1,w_2)$
in the terms of $\xi.$  Let us start with the deterministic
admissible shifts. Denote by $i$ the canonical embedding of $H$ into
$C^2,$ i.e.
$$
i(\vf)(t)=(\int^t_0\vf_1ds,\int^t_0\vf_2ds).
$$

{\bf Lemma 4.1.} {\it Let the operator norm $\|V\|<1,$  then for
every $h\in H$  $i(h)$ is admissible shift for $\mu_{w_1,w_2}$  and
the corresponding density has the form
$$
p_h(\xi)=\exp\{(S^{-1}h,\xi)-\frac{1}{2}(S^{-1}h,h)\}. \eqno(2.2)
$$
}

{\bf Remark.}   Due to the condition $\|V\|<1$ the operator $S^{-1}$
is bounded on $H$  and can be written in the form
$$
S^{-1}\vf=\vf+(Q_{11}\vf_1+Q_{12}\vf_2,
Q_{21}\vf_1+Q_{22}\vf_2)=\vf+Q\vf,
$$
where $\|Q\|<1.$

{\it Proof.}  Note, that by the definition the operator $S$ is
nonnegative. Define
$$
{\xi}'=S^{-\frac{1}{2}}\xi.
$$
Then the shift of the distribution of $(w_1,w_2)$ on the vector
$i(h)$  is related to the shift of ${\xi}'$  on the vector
$S^{-\frac{1}{2}}h.$   Now the statement of the lemma follows from
the well-known formula for the density in the terms of ${\xi}'$
$$
p({\xi}')=\exp\{({\xi}',
S^{-\frac{1}{2}}h)-\frac{1}{2}(S^{-\frac{1}{2}}h,
S^{-\frac{1}{2}}h)\}
$$
if we rewrite it in the terms of $\xi.$

The lemma is proved.

{\bf Remark.}  Formula (2.2) can be rewritten in the terms of
$w_1,w_2.$ Really, by the definition
$$
(S^{-1}h,\xi)-\frac{1}{2}(S^{-1}h,h)=
$$
$$
= \int^1_0(S^{-1}h)_1dw_1+ \int^1_0(S^{-1}h)_2dw_2- \frac{1}{2}
\int^1_0(S^{-1}h)_1h_1ds- \frac{1}{2} \int^1_0(S^{-1}h)_2h_2ds=
\eqno(2.3)
$$
$$
= \int^1_0(h_1+Q_{11}h_1+Q_{12}h_2)dw_1+
\int^1_0(h_2+Q_{21}h_1+Q_{22}h_2)dw_2-
$$
$$
-\frac{1}{2} \int^1_0(h_1+Q_{11}h_1+Q_{12}h_2)h_1ds- \frac{1}{2}
\int^1_0(h_2+Q_{21}h_1+Q_{22}h_2)h_2ds.
$$

Using the same method one can find the density of $\mu_{x_1,x_2}$
with respect  $\mu_{w_1,w_2}.$  Firstly define the stochastic
derivatives of the functionals from $w_1,w_2$  with respect to $\xi$
and the extended stochastic integral in the terms of $\xi.$  Let
$\vf$  be a differentiable bounded function on $C\oplus C.$   Define
the stochastic derivative of the random variable $\vf(w_1,w_2)$  by
the formula
$$
D\vf(w_1,w_2):=i^*\nabla\vf(w_1,w_2).
$$
By this definition for every $t\in[0;1]$
$$
Dw_1(t)=(\1_{[0;t]},0),
$$
$$
Dw_2(t)=(0, \1_{[0;t]}).
$$
Note, that $\vf(w_1,w_2)$   can be regarded as a functional from the
generalized random element ${\xi}'$  which was introduced in the
proof of the lemma 4.1. Since ${\xi}'$  has an identity covariation
operator the stochastic derivatives and extended stochastic integral
for the functionals from ${\xi}'$  are connected by the usual
relation. Now we will define the extended stochastic integral with
respect to $\xi.$  It can be done in the following way. Consider the
Gaussian random functional on $H$  of the kind
$$
J(\vf)=(\vf, \xi).
$$
Then, in the terms of ${\xi}'$  $J$  can be rewritten as
$$
J(\vf)=(S^{\frac{1}{2}}\vf, {\xi}').
$$
So, the action of $J$  on the random element $x$  in $H$  via the
definition 1.3 has the form
$$
J(x)=I(S^{\frac{1}{2}}x). \eqno(2.4)
$$
Here $I$ is the extended stochastic integral with respect to
${\xi}'.$ Note also, that for the stochastic derivatives with
respect to ${\xi}'$  and ${\xi}$ we have the obvious relation
$$
D_{{\xi}'}\alpha=S^{-\frac{1}{2}}D_\xi\alpha.
$$
Hence on the domain of definition
$$
E(D_\xi\alpha, x)=E(S^{-\frac{1}{2}}D_\xi\alpha, S^{\frac{1}{2}}x)=
$$
$$
= E(D_{{\xi}'}\alpha, S^{\frac{1}{2}}x) =E\alpha\cdot
I(S^{\frac{1}{2}}x)= E\alpha J(x). \eqno(2.5)
$$
Thus the relation between the stochastic derivative and extended
stochastic integral with respect to $\xi$ is the same as for
${\xi}'.$ Now let us turn to the nonadapted shifts of the
distribution of $(w_1,w_2).$ Consider the pair of random processes
$x_1,x_2$  which are defined by the equations (2.1).

 The next lemma is standard.

{\bf Lemma 2.2.} {\it Let the functions $a_1,a_2$  be continuously
differentiable and have bounded derivatives. Then

1) for every $t\in[0;1]$  the random variables $x_1(t), x_2(t)$
have the stochastic derivatives $Dx_1(t), Dx_2(t),$

2) the random element $(a_1(x_1(\cdot)), a_2(x_1(\cdot)))$ in $H$
has the stochastic derivative, and
$$
D(a_1(x_1(s)), a_2(x_1(s)))(t)=
$$
$$
= (a'_1(x_1(s))Dx_1(s)(t), a'_2(x_1(s))Dx_1(s)(t)),
$$

3) the stochastic derivative of $x_1$  with respect to $w_1$  (i.e.
the first coordinate of $Dx_1$)  satisfies the equation
$$
\begin{array}{l}
D_1x_1(s)(t)=1+\int^s_ta'_1(x_1(r))D_1x_1(r)(t)dr, 0\leq t\leq s\leq1,\\
D_1x_1(s)(t)=0, \ \ t>s,
\end{array}
$$
and
$$
Dx_1(s)(t)=(D_1x_1(s)(t),0).
$$
}

It follows from the lemma 2.2, that $\|D(a_1(x_1(\cdot)),
a_2(x_1(\cdot)))\|_H$ can be made small if we take $a_1'$ and $a_2'$
small enough. Since the operator $S^{-\frac{1}{2}}$  is bounded in
$H,$  then due to the theorem 3.2.2 from [15]  the distribution of
$(x_1,x_2)$  is absolutely continuous with respect to the
distribution $(w_1,w_2)$   for sufficiently small $a_1',$ $a_2'.$
The corresponding density will be denoted by $p.$ Accordingly to
[18]  $p$  has the form
$$
p=\zeta\cdot\exp\{I(S^{-\frac{1}{2}}h)-\frac{1}{2}(S^{-1}h,h)\},
\eqno(2.6)
$$
where
$$
h(t)=(a_1(w_1(t)), a_2(w_1(t))),
$$
and $\zeta$  is the corresponding Carleman--Fredgolm determinant.
Due to (2.4) (2.6) can be rewritten as
$$
p=\zeta\exp\{J(S^{-1}h)-\frac{1}{2}(S^{-1}h, h)\}. \eqno(2.7)
$$
This expression allows us to conclude, that up to the term $\zeta$
$p$  has the stochastic derivative. We will suppose, that this is so
in the next section, where the formulas for the conditional
expectation and extended stochastic integral will be obtained in
non-Gaussian case. \vskip15pt

{\bf3. Conditional expectation.}  For the processes $(x_1,x_2)$ from
(2.1) let us search for the conditional distribution of $x_1(t)$
under fixed $\{x_2(s); s\in[0;1]\}.$   Firstly note, that under our
conditions the distribution of $(x_1,x_2)$  is absolutely continuous
with respect to the distribution of $(w_1,w_2)$  and, consequently,
the distribution  of $x_2$  is absolutely continuous with respect to
the distribution of $w_2.$

Denote for a moment by $\mu$  the distribution of the pair
$(w_1,w_2)$  on $C([0;1])\oplus C([0;1])$   and by $\mu_1$  and
$\mu_2$  the distributions of $w_1$  and $w_2$  (surely these are a
Wiener measures but on the different  copies of $C([0;1]).$ It
follows from the general theory of integration   that the measure
$\mu$   can be desintegrated with respect to $\mu_2,$  i.e.
$$
\mu(\Delta)=\int_{C([0;1])}\nu(u,\Delta_u)\mu_2(du),
$$
for arbitrary Borel $\Delta$  in $C([0;1])\oplus C([0;1]).$  Here
$\nu$  is a measurable family of the probability measures and
$\Delta_u=\{v\in C([0;1]): (v,u)\in \Delta\}.$

Define for the measurable bounded function $\vf: C([0;1])\to\mbR$
the function $\psi$  by the rule
$$
C([0;1])\ni u\mapsto\psi(u)=\int_{C([0;1])}\vf(v)p(v,u)\nu(u,dv).
\eqno(3.1)
$$
$$
\cdot (\int_{C([0;1])}p(v,u)\nu(u, dv))^{-1}.
$$

The following variant of the Bayes formula holds.

{\bf Lemma 3.1.}
$$
E(\vf(x_1)/x_2)=\psi(x_2).
$$

{\it Proof.}   Note firstly, that $\psi(x_2)$  is correctly defined
because the function $\psi$  is defined up to the set of Wiener
measure zero, and the distribution of $x_2$ is equivalent to this
measure. Now for arbitrary bounded and  measurable function $\gamma:
C([0;1])\to\mbR$
$$
E\vf(x_1)\gamma(x_2)=E\vf(w_1)\gamma(w_2)p(w_1,w_2)=
$$
$$
=E\gamma(w_2)\cdot E(\vf(w_1)p(w_1,w_2)/w_2)=
$$
$$
 =E\gamma(w_2)\cdot \frac
{E(\vf(w_1)p(w_1,w_2)/w_2)}{E(p(w_1,w_2)/w_2)}\cdot
E(p(w_1,w_2)/w_2)=
$$
$$
= E\gamma(w_2)\psi(w_2)\cdot p(w_1,w_2)=E\gamma(x_2)\psi(x_2).
$$
This finishes the proof.

For arbitrary $t\in[0;1]$  denote by $\pi_t$   the random measure on
$\mbR$   whose pairing with the bounded measurable function $f$  is
defined by the formula
$$
\int_{\mbR} f(r)\pi_t(dr)=E(f(w_1(t))\cdot p(w_1,w_2)/w_2).
$$

In view of the previous lemma it is enough to get the equation for
$\pi_t.$ The next lemma contains the necessary facts from the theory
of extended stochastic integral.

{\bf Lemma 3.2.} {\it Let $H$  be the separable Hilbert space, $\xi$
be a generalized Gaussian random element in $H$ with zero mean and
identity covariation. Suppose, that the random element $x$  in $H$
has two stochastic derivatives and let $I$ and $D$   be the symbols
of the extended stochastic integral and stochastic derivative
correspondingly.  Then for arbitrary $h\in H$  and stochastically
differentiable bounded random variable $\alpha$   the following
formulas hold

1) $\alpha I(x)=I(\alpha x)+(x, D\alpha),$

2) $(DI(x), h)=(x,h)+I((Dx, h)).$ }

{\it Proof.}   The first statement is the well-known relation [12].
Let us check 2).  Use the integration by part formula. Consider the
random variable $\beta$  which is twice stochastically
differentiable. Then, using 1),
$$
E(DI(x),h)\beta =E(DI(x), \beta h)=
$$
$$
= EI(x)I(\beta h) =EI(x)(\beta I(h)-(D\beta, h))=
$$
$$
=E(x, D(\beta I(h)-(D\beta,h)))=
$$
$$
= E[(x, \beta h)+(x, D\beta)I(h)-(x, (D^2\beta, h))]=
$$
$$
= E(x,h)\beta +E(x, I(D\beta h))=
$$
$$
=E\beta((x,h)+I((Dx,h))).
$$
Lemma is proved.

{\bf Remark.}   Note, that the statement of the lemma remains to be
true in the case, when $\alpha$  is not bounded but all terms are
well-defined. Also, due to the formula (4.5) the lemma holds in the
case, when the initial generalized Gaussian random element has not
identity covariation.

Now let us turn  to our filtration problem.

Take the function $f\in C^2_0(\mbR).$  Then for arbitrary $r\in\mbR$
from the Ito formula
$$
f(r+w_1(t))p(w_1,w_2)=
$$
$$
=f(r)p(w_1,w_2)+\int^t_0f'(r+w_1(s))dw_1(s)p(w_1,w_2)+
$$
$$
+\frac{1}{2}\int^t_0f''(r+w_1(s))p(w_1,w_2)ds.
$$
Consider the second summand. It contains the Ito integral which
considers with the extended stochastic integral as it was mentioned
before. So we can apply the  formula 1) from the lemma 3.2:
$$
p(w_1,w_2)\int^t_0f'(r+w_1(s))dw_1(s)=
$$
$$
=\int^t_0f'(r+w_1(s))p(w_1,w_2)dw_1(s)+
$$
$$
+\int^t_0f'(r+w_1(s))(SDp(w_1,w_2))_1(s)ds.
$$
Here the index 1 symbolize the first coordinate of correspondent
element from $H.$  Now note,  that the conditional expectation with
respect to $w_2$  is an operator of the second quantization. So, if
we denote by
$$
\gamma_1(t)=E(w_1(t)/w_2),
$$
then, due to the theorem 1.1
$$
E(f(r+w_1(t))p(w_1,w_2)/w_2)=
$$
$$
=
E(f(r)p(w_1,w_2)/w_2)+\int^t_0E(f'(r+w_1(s))p(w_1,w_2)/w_2)d\gamma(t)+
$$
$$
+\frac{1}{2} \int^t_0E(f''(r+w_1(s))p(w_1,w_2)/w_2)ds+ \eqno(3.2)
$$
$$
+\int^t_0E(f'(r+w_1(s)(SDp(w_1,w_2))_1)(s)/w_2)ds,
$$
where the integral with respect to $\gamma$ is an extended
stochastic integral.
 In order to get the stochastic differentiability of $p$  let us consider
the case when Carleman-Fredholm determinant $\zeta$  is equal to
one. Denote by $P_t$  the orthogonal projector in $L_2([0;1])\oplus
L_2([0;1])$ on the subspace $L_2([0;t])\oplus L_2([0;t]).$

{\bf Lemma 3.3.} {\it Suppose, that the operator $S$  has the
property
$$
\forall t\in[0;1]: \ \ \ \  P_tS=P_tSP_t.
$$
Then $\zeta=1.$ }

{\it Proof.} The value $\zeta$ is the Carleman-Fredholm determinant
of the operator $SDh,$   where
$$
h=(a_1(w_1(\cdot)), a_2(w_1(\cdot))).
$$
Now
$$
Dh(t,s)=
\begin{pmatrix}
a'_1(w_1(t))\1_{[0;t]}(s)&0\\
\\
a'_2(w_1(t))\1_{[0;t]}(s)&0
\end{pmatrix}
\eqno(3.3)
$$

In order to prove that
$$
\det\nolimits_2(Id+SDh)=1
$$
we will use the theorem 3.6.1 from [15].  Due to this theorem it is
enough to check, that the operator $SDh$   is quasi-nilpotent, i.e.,
that
$$
\lim_{n\to\infty}\|(SDh)^n\|^{\frac{1}{n}}=0. \eqno(3.4)
$$
It follows from representation (3.3), that\newline $\forall \vf\in H
\ \forall t\in[0;1]$
$$
\|P_tDh\vf\|^2\leq c\int^t_0\|P_s\vf\|^2ds,
$$
where $c$  depends on the $\sup_\mbR(|a'_1|+|a'_2|).$

Consequently
$$
\|(SDh)^n(\vf)\|^2\leq\|S\|^2\cdot c\cdot
\int^1_0\|P_{t_1}(SDh)^{n-1}(\vf) \|^2dt_1\leq
$$
$$
\leq \|S\|^2\cdot c\cdot\int^1_0 \|P_{t_1}S
P_{t_1}P_{t_1}Dh(SDh)^{n-2}(\vf)\|^2dt_1\leqq
$$
$$
\leq \|S\|^4\cdot c^2\cdot\int^1_0\int^{t_1}_0
\|P_{t_2}(SDh)^{n-2}(\vf)\|^2dt_2dt_1\leq\ldots\leq
$$
$$
\leq \|S\|^{2n}\cdot
c^n\cdot\int^1_0\int^{t_1}_0\ldots\int^{t_{n-1}}_0
\|P_{t_n}\vf\|^2dt_n\ldots dt_1\leq
$$
$$
\leq\|\vf\|^2\frac{\|S\|^{2n}c^n}{n!}.
$$
This means, that (3.4) holds and $\zeta=1.$

Lemma is proved.

Now one can conclude that $p$  has the stochastic derivative and
(3.1)  is correct. The further concretization of (3.2)  can be
possible due to the special form of $p.$   As a  consequence of
(3.2) and (3.4)   we have the following theorem.

{\bf Theorem 3.1.}  {\sl Suppose that the coefficients $a_1, a_2$
and the operator $V$ satisfy the conditions of the lemma 3.3. Then
the random function
$$
U(r,t)=E(f(r+w_1(t))p(w_1,w_2)/w_2)
$$
satisfies relation
$$
dU(r,t)=\frac{1}{2} \frac{\pt^2}{\pt r^2}U(r,t)dt+ \eqno(5.5)
$$
$$
+\frac{\pt}{\pt r} U(r,t)\gamma(dt)+Ef'(r+w_1(t))
(SDp(w_1,w_2))_1(t)dt.
$$
}

In some particular case the last term can be written in a simple
form. For example, when  $a_2=0,$ then (3.5)  transforms into
$$
dU(r,t)=\frac{1}{2} \frac{\pt^2}{\pt r^2} U(r,t)dt+\frac{\pt}{\pt
r}U(r,t)\gamma(dt)+
$$
$$
+a_1(r)\frac{\pt}{\pt r}U(r,t)dt.
$$

\vskip 1cm \centerline{REFERENCES}

\begin{enumerate}
\item
It\^o, Kiyosi; McKean, Henry P., Jr. Diffusion processes and their
sample paths. Second printing, corrected. Die Grundlehren der
mathematischen Wissenschaften, Band 125. Springer-Verlag,
Berlin-New York, 1974.
\item
Simon Barry. The $P(\phi)_2$ euclidian (quantum) field theory.
Princeton University Press, 1974.
\item
 Dorogovtsev A.A. Stochastic analysis and random maps in Hilbert
space. VSP Utrecht, The Netherlands, Tokyo, Japan. - 1994. - 110
p.
\item
Dorogovtsev A.A. Anticipating equations and filtration problem.
Theory of stochastic processes, 1997. V.3(19), issue 1 - 2. P. 154
- 163.
\item
Dorogovtsev A.A. Conditional measures for diffusion processes and
anticipating stochastic equations. Theory of stochastic processes,
1998. V.4(20). P. 17 - 24.
\item
 Skorokhod A.V. One generalization of the stochastic integral.
Probability theory and its applications. 1975, V. 20, \#2. P.223 -
237.
\item
 Liptser R. Sh., Shyriaev A.N. Statistics of random processes.
Moskow, Nauka, 1974. 696 p.
\item
Billingsley Patrick. Convergence of probability measures. John
Wiley and Sons, Inc., New York, London, Sydney, Toronto. 1968.
\item
Mandelbrot, B.B., van Ness, J. Fractional Brownian motion,
fractional noises and applications, SIAM Rev., 10:422-437. 1968.
\item
 Tindel, S.; Tudor, C.A.; Viens, F. Stochastic evolution
equations with fractional Brownian motion. Probab. Theory Relat.
Fields 127, No. 2, 186-204.
\item
Dorogovtsev A.A. An action of Gaussian strong random operator on
random elements. Probability theory and its  applications. - 1986.
- V.31,N 4. - P. 811-814.
\item Malliavin, Paul Stochastic analysis. 1997, Text.Mo\-no\-graph,
Grundlehren der Mathematischen Wissenschaften. 313., Berlin:
Springer.
\item
 Nualart, David The Malliavin calculus and related topics. 1995,
Text.Mono\-graph, Probability and Its Applications., New York, NY:
Springer-Verlag.
\item
 Dorogovtsev A.A. Stochastic integration and one class of
Gaussian stochastic processes. Ukr. math. journal, 1998. V.50,\#4.
P.495 - 505.
\item
A.Suleyman Ustunel, Moshe Zakai.  Transformation of Measure on
Wiener Space. -- Springer, 2000. -- 298 p.
\end{enumerate}

\end{document}